# AUTOMORPHISM GROUPS FOR SEMIDIRECT PRODUCTS OF CYCLIC GROUPS

JASHA SOMMER-SIMPSON

ABSTRACT. Every semidirect product of groups $K \rtimes H$ has size $|K| \cdot |H|$, yet the size of such a group's automorphism group varies with the chosen action of $H$ on $K$. This paper will explore groups of the form $\mathrm{Aut}(K \rtimes H)$, considering especially the case where $H$ and $K$ are cyclic. Only finite groups will be considered.

## Contents

## 1. Introduction

There exist some finite groups that are isomorphic to their own automorphism groups, e.g. $D_6$. The structure of $\mathrm{Aut}(G)$ is often hard to compute, and therefore it is difficult to decide in general when $G \cong \mathrm{Aut}(G)$.

In some cases the order of $\mathrm{Aut}(G)$ is easy to compute even while the structure of $\mathrm{Aut}(G)$ remains enigmatic. For example, the structure of $\mathrm{Aut}(\mathbb{Z}_n \times \mathbb{Z}_2)$ is quite complicated when $n \equiv 2 \pmod{4}$, but the orders of these automorphism groups are easily calculable (see Table 1). This paper will explore the orders and structures of automorphism groups for semidirect products, focusing specifically on the case $\mathrm{Aut}(\mathbb{Z}_m \rtimes \mathbb{Z}_n)$.

| The value of $n$ modulo 4 | The order of $\mathrm{Aut}(\mathbb{Z}_n \times \mathbb{Z}_2)$ |
|---|---|
| $n \equiv 0 \pmod{4}$ | $4\phi(n)$ |
| $n \equiv 1 \pmod{4}$ | $\phi(n)$ |
| $n \equiv 2 \pmod{4}$ | $6\phi(n)$ |
| $n \equiv 3 \pmod{4}$ | $\phi(n)$ |

TABLE 1. $\phi$ denotes the Euler totient function.

Section 2 will introduce the notation employed via an informal discussion of ideas from group theory, and section 3 will provide a number of well-known definitions and theorems concerning sequences and abelian groups.

In section 4 we will compute the order of $\mathrm{Aut}(G)$ for a several abelian groups $G$, e.g. the arbitrarily large direct product $\mathbb{Z}_p \times \mathbb{Z}_p \times \cdots \times \mathbb{Z}_p$ where $p$ is prime. The results from Table 1 will be covered in depth. Section 5 gives general results concerning $\mathrm{Aut}(H \times K)$, while section 6 uses information about the action $\psi \in \mathrm{hom}\big(H, \mathrm{Aut}(K)\big)$ to analyze $\mathrm{Aut}(K \rtimes_\psi H)$.

In section 7 we ask, "When is $G$ isomorphic to its own automorphism group?" While no general solution is provided, we are able to answer this question in the dihedral case. A general form for the automorphisms and automorphism groups of dihedral groups will be provided. The final part of this paper (section 8) will describe $\mathbb{Z}_8 \rtimes \mathbb{Z}_2$ for each action of $\mathbb{Z}_2$ on $\mathbb{Z}_8$, and will compute explicitly the automorphism groups for these semidirect product groups. This should serve as an interesting example because several distinct $\mathbb{Z}_8 \rtimes \mathbb{Z}_2$ groups have isomorphic automorphism groups.

While it is difficult in general to find the structure of $\mathrm{Aut}(\mathbb{Z}_m \rtimes \mathbb{Z}_n)$, or indeed even to visualize the multiplication in $\mathbb{Z}_m \rtimes \mathbb{Z}_n$ for various actions of $\mathbb{Z}_n$ on $\mathbb{Z}_m$, I hope that this paper may serve to shed light on the problem via examples with groups of small order.

## 2. Notation

Capital letters represent groups, and Greek letters will in general be group homomorphisms. However, the letter $\phi$ denotes the Euler totient function, to be defined after this paragraph. If $\psi$ is a group homomorphism from $A$ to $B$ and the elements of $B$ are themselves group homomorphisms, $\psi_a$ will sometimes be used to represent the image in $B$ of $a \in A$. Given a natural number $n$, we will write $[1, n]$ for the inclusive set ot integers ranging from 1 to $n$.

For any element $g \in G$, let $\langle g \rangle$ be the cyclic subgroup of $G$ generated by $g$, and let $|g|$ be the order of $g$. The center of $G$, given by $Z(G) = \{x \in G \mid gx = xg \; \forall g \in G\}$, is the subgroup containing elements that commute with every other element in $G$. The automorphism group of $G$, written $\mathrm{Aut}(G)$, has elements given by the isomorphisms from $G$ to $G$. Multiplication in this group is defined by composition of functions.

In general, we will write $e$ for a group's identity element, and we will use $\{e\}$ to represent the trivial group of order 1. However, we will sometimes use $\mathrm{id}_G$ to refer to the identity element in $\mathrm{Aut}(G)$.

Given $n \in \mathbb{N}$, the totatives of $n$ are the positive integers $x$ such that $x \leq n$ and $\gcd(x, n) = 1$. The Euler totient function $\phi$ takes input from the natural numbers and is defined by $\phi(n) = \big|\{x \in \mathbb{N} \mid 1 \leq x \leq n \text{ and } \gcd(x, n) = 1\}\big|$. Specifically, $\phi(n)$ counts the totatives of $n$. Note that $\phi(1) = 1$, as $\gcd(1, 1) = 1$. This function is particularly useful in considering the automorphism groups for cyclic groups: given a generator $r$ of $\mathbb{Z}_n$, we know that $\langle r^i \rangle \cong \mathbb{Z}_n$ if and only if $\gcd(i, n) = 1$. Therefore there are as many generators for $\mathbb{Z}_n$ as there are relative primes to $n$ in $[1, n]$. Hence the order of $\mathrm{Aut}(\mathbb{Z}_n)$ equals $\phi(n)$.

Let $\mathrm{hom}(X, Y)$ be the set of all group homomorphisms from $X$ to $Y$. For an element $\mu$ in $\mathrm{hom}(X, Y)$, we say that $\mathrm{im}(\mu)$ is the subgroup of $Y$ given by the image of $X$ through $\mu$, i.e. $\{\mu(x) \in Y \mid x \in X\}$. Furthermore, we define $\ker(\mu)$ to be the kernel of $\mu$, given by $\{x \in X \mid \mu(x) = e \in Y\}$. This is subgroup of $X$ whose elements map to the identity element in $Y$.

Given a homomorphism $\psi : H \to \operatorname{Aut}(K)$, we write $K \rtimes_\psi H$ for the semidirect product of $K$ by $H$. This is a group whose elements look like those of $K \times H$, although the multiplication in $K \rtimes_\psi H$ is defined differently: the product of $(k, h)$ and $(k', h')$ is given by $(k\psi_h(k'), hh')$ where $\psi_h$ is the image of $h$ through $\psi$. Relaxing the requirement for parentheses, one could write $khk'h' = k\psi_h(k')hh'$.

## 3. Sequences and Abelian Groups

While the topics discussed here may not be as basic as those from section 2, readers familiar with group theory will already be acquainted with much this section's material. These definitions and theorems are used implicitly throughout the paper.

### 3.1. Short Exact Sequences.

Given group homomorphisms $\iota : K \to G$ and $\tau : G \to H$, the sequence of groups $K \xrightarrow{\iota} G \xrightarrow{\tau} H$ is said to be <u>exact at $G$</u> if $\operatorname{im}(\iota) = \ker(\tau)$. The map $\iota$ is injective if and only if the sequence $\{e\} \to K \xrightarrow{\iota} G$ is exact at $K$, while $\tau$ is surjective if and only if we have exactness at $H$ in $G \xrightarrow{\tau} H \to \{e\}$. Given $\iota$ injective, $\tau$ surjective, and $\operatorname{im}(\iota) = \ker(\tau)$, we say that $\{e\} \to K \xrightarrow{\iota} G \xrightarrow{\tau} H \to \{e\}$ is a <u>short exact sequence</u>. This means that we have exactness at $K$, at $G$ and at $H$.

Suppose that $\{e\} \to K \xrightarrow{\iota} G \xrightarrow{\tau} H \to \{e\}$ is a short exact sequence. Because $\iota$ is injective we can say that $K$ is a subgroup of $G$, equal to to $\operatorname{im}(\iota)$ and $\ker(\tau)$. Applying the First Isomorphism Theorem, we see that $\operatorname{im}(\tau) \cong G/\ker(\tau)$, or equivalently $H \cong G/K$. If there exists a homomorphism $\rho : H \to G$ satisfying $(\tau \circ \rho) = \operatorname{id}_H$, then our short exact sequence is said to <u>split</u>. This map $\rho$ is called a <u>section</u>. Notice that $\rho$ is necessarily injective, so we can say that $H$ is a subgoup of $G$ equal to $\operatorname{im}(\rho)$.

**Theorem 3.1.** *The short exact sequence $\{e\} \to K \xrightarrow{\iota} G \xrightarrow{\tau} H \to \{e\}$ splits if and only if $G \cong K \rtimes_\psi H$ for some $\psi \in \hom\big(H, Aut(K)\big)$.*

*Proof.* First suppose that we are given an isomorphism $\Gamma : K \rtimes_\psi H \to G$. We will construct a short exact sequence $\{e\} \to K \xrightarrow{\iota} G \xrightarrow{\tau} H \to \{e\}$ and a section $\rho : H \to G$. We can define injective homomorphisms $\iota : K \to G$ and $\tau : H \to G$ by $k \mapsto \Gamma(k, e)$ and $h \mapsto \Gamma(e, h)$, respectively.

As $\Gamma$ is an isomorphism, we can say that every element of $G$ has the form $\Gamma(k, h)$ for some $(k, h)$ in $K \rtimes_\psi H$. To define a surjective homomorphism $\tau : G \to H$, let $\tau\big(\Gamma(k, h)\big) = h$ for every element $\Gamma(k, h)$ in $G$. We need now to verify that $\operatorname{im}(\iota) = \ker(\tau)$ and that $(\tau \circ \rho) = \operatorname{id}_H$. Indeed, $(\tau \circ \rho)(h) = \tau\big(\Gamma(e, h)\big) = h \ \forall \ h \in H$, and $\operatorname{im}(\iota) = \{\Gamma(k, e) \in G\} = \{\Gamma(k, h) \in G \mid \tau\big(\Gamma(k, h)\big) = e\} = \ker(\tau)$.

Now suppose that we are given a short exact sequence $\{e\} \to K \xrightarrow{\iota} G \xrightarrow{\tau} H \to \{e\}$ and a section $\rho$. We must find a homomorphism $\psi : H \to \operatorname{Aut}(K)$ that allows for the construction of an isomorphism $\Theta : G \to K \rtimes_\psi H$.

Let $\pi : \operatorname{im}(\iota) \to K$ be the isomorphism such that $\pi \circ \iota = \operatorname{id}_K$ and $\iota \circ \pi = \operatorname{id}_{\operatorname{im}(\iota)}$. Because $\operatorname{im}(\iota)$ is normal in $G$, the element $\tau(h)\iota(k)\tau(h^{-1})$ is in $G$ for all $k \in K$ and $h \in H$. Thus for every $h$ in $H$ we can define an isomorphism $\psi_h : K \to K$ by $k \mapsto \pi\big(\tau(h)\iota(k)\tau(h^{-1})\big)$. Moreover, the map $\psi : H \to \operatorname{Aut}(K)$ given by $h \mapsto \psi_h$ is

homomorphic because

$$(\psi_{h_1} \circ \psi_{h_2})(k) = \psi_{h_1}\Big(\pi\big(\tau(h)\iota(k)\tau(h^{-1})\big)\Big) = \pi\bigg(\tau(h_1)\iota\Big(\pi\big(\tau(h_2)\iota(k)\tau(h_2^{-1})\big)\Big)\tau(h_1^{-1})\bigg)$$

$$= \pi\big(\tau(h_1)\tau(h_2)\iota(k)\tau(h_2^{-1})\tau(h_1^{-1})\big) = \pi\big(\tau(h_1 h_2)\iota(k)\tau(h_2^{-1}h_1^{-1})\big) = \psi_{h_1 h_2}(k)$$

for all $h_1, h_2 \in H$ and $k \in K$.

We know that $|G| = |K||H|$ and that every right coset of $\mathrm{im}(\iota)$ in $G$ can be written as $K\rho(h)$ for some $h \in H$. Therefore every element $x$ in $G$ can be expressed as $\iota(k)\rho(h)$ for some unique $k \in K$ and $h \in H$. Specifically, $h = \tau(x)$ and $k = \pi\big(x\rho(h^{-1})\big)$.

To see that the function $\Theta : G \to K \rtimes_\psi H$ given by $\iota(k)\rho(h) \mapsto kh$ is an isomorphism, we need only to check that $\Theta(xy) = \Theta(x)\Theta(y)$ for all $x, y \in G$. Without loss of generality, write $x = \iota(k)\rho(h)$ and $y = \iota(k')\rho(h')$.

$$\Theta\big(\iota(k)\rho(h)\iota(k')\rho(h')\big) = \Theta\big(\iota(k)\rho(h)\iota(k')\rho(h^{-1})\rho(h)\rho(h')\big) = \Theta\Big(\iota(k)\iota\big(\psi_h(k')\big)\rho(hh')\Big)$$

$$= k\psi_h(k')hh' = khk'h' = \Theta\big(\iota(k)\rho(h)\big)\Theta\big(\iota(k')\rho(h')\big)$$

$\square$

**Theorem 3.2.** *If $H$ is a subgroup of $G$ such that $\frac{|G|}{|H|} = 2$ then $H \triangleleft G$.*

*Proof.* Suppose that the subgroup $H$ contains half the elements from $G$. For all $x \in H$ we have $xH = H = Hx$. When $x \in G$ but $x \notin H$, the cosets $xH$ and $Hx$ must contain precisely the same elements: those from $G$ that are not in $H$. Therefore $xH = Hx$ for all $x \in G$. $\square$

3.2. **Pertaining to Cyclic Groups.**
The theorem below will be useful in considering direct products of cyclic groups with relatively prime order.

**Theorem 3.3** (Chinese Remainder Theorem)**.** *Suppose that $n = p_1^{k_1} p_2^{k_2} \cdots p_w^{k_w}$ for distinct prime numbers $p_1, p_2, \ldots, p_w$ and for natural numbers $k_1, k_2, \ldots, k_w$. Then $\mathbb{Z}_n$ is isomorphic to $\mathbb{Z}_{p_1^{k_1}} \times \mathbb{Z}_{p_2^{k_2}} \times \cdots \times \mathbb{Z}_{p_w^{k_w}}$.*

*Proof.* Write $G = \mathbb{Z}_{p_1^{k_1}} \times \mathbb{Z}_{p_2^{k_2}} \times \cdots \times \mathbb{Z}_{p_w^{k_w}}$. As $G$ has order $n = p_1^{k_1} p_2^{k_2} \cdots p_w^{k_w}$, we can prove the isomorphism $G \cong \mathbb{Z}_n$ by finding an element in $G$ that has order $n$.

Let $\{g_1, g_2, \ldots, g_w\}$ be a set of generators for $G$ such that each $g_t$ generates the subgroup $\mathbb{Z}_{p_t^{k_t}}$ of $G$. To calculate the order of the element $g_1 g_2 \cdots g_w$ in $G$, we must find the smallest positive integer $s$ satisfying $(g_1 g_2 \cdots g_w)^s = e$. But $(g_1 g_2 \cdots g_w)^s = e$ if and only if $g_x^s = e$ for every $x \in [1, w]$.

Thus we have $(g_1 g_2 \cdots g_w)^s = e$ if and only if $s \equiv 0 \pmod{p_x^{k_x}}$ for every $x \in [1, w]$. As the numbers $p_1^{k_1}, p_2^{k_2}, \ldots, p_w^{k_w}$ are pairwise relatively prime, the smallest positive integer $s$ satisfying

$$s \equiv 0 \pmod{p_x^{k_x}} \text{ for every } x \in [1, w]$$

is $n = p_1^{k_1} p_2^{k_2} \cdots p_w^{k_w}$, and therefore the order of $g_1 g_2 \cdots g_w$ is equal to $n$. $\square$

**Theorem 3.4** (Classification of Finite Abelian Groups)**.** *A finite group $G$ is abelian if and only if $G$ is isomorphic to the direct product of finitely many cyclic groups.*

**Theorem 3.5.** *Given two elements $x$ and $y$ in an abelian group $G$, the order of $xy$ is $lcm(|x|, |y|)$.*

*Proof.* The order of $xy$ in an abelian group $G$ is the smallest positive integer $s$ satisfying $(xy)^s = e$. The abelian nature of $G$ gives that $(xy)^s = x^s y^s$, which means we must find the smallest $s$ that simultaneously satisfies $x^s = e$ and $y^s = e$. Therefore $(xy)^s = e$ if and only if both $|x|$ and $|y|$ divide $s$; the smallest such $s$ is $lcm(|x|, |y|)$. □

**Theorem 3.6** (Automorphism Groups of Cyclic Groups). $Aut(\mathbb{Z}_n) \cong (\mathbb{Z}_n)^\times$

*Proof.* Suppose that $r$ generates $\mathbb{Z}_n$. The multiplicitive group of invertible integers modluo $n$ has elements of the form $r^i$ for $i$ relatively prime to $n$. Multiplication in this group is given by $r^i r^j = r^{ij}$.

An element $r^i$ of $\mathbb{Z}_n$ generates $\mathbb{Z}_n$ if and only if $\gcd(i, n) = 1$. Therefore every element of $Aut(\mathbb{Z}_n)$ must send $r$ to $r^i$ for some totative $i$ of $n$.

Consider the map $\nu : (\mathbb{Z}_n)^\times \to Aut(\mathbb{Z}_n)$ sending each $r^i$ in $(\mathbb{Z}_n)^\times$ to the automorphism $\nu_{r^i} : \mathbb{Z}_n \to \mathbb{Z}_n$ defined by $r \mapsto r^i$. First, notice that $\nu$ is bijective. Then we can observe that $\nu$ is homomorphic because
$$\nu_{r^i r^j}(r) = \nu_{r^{ij}}(r) = r^{ij} = (r^j)^i = \nu_{r^i}(r^j) = (\nu_{r^i} \circ \nu_{r^j})(r).$$
□

**Theorem 3.7.** *When $n$ is odd, $\phi(n) = \phi(2n)$.*

**Theorem 3.8.** *For $p$ prime, $Aut(\mathbb{Z}_p) \cong \mathbb{Z}_{(p-1)}$.*

Given that this paper concerns automorphisms, the following definition will be useful in understanding the groups to which our automorphisms pertain.

**Definition 3.9.** Let $C < G$. $C$ is called a characteristic subgroup of $G$ if every map in $Aut(G)$ sends the elements of $C$ amongst themselves.

Notice that if $C$ is characteristic in $G$ then $C$ is normal in $G$.

## 4. Aut Groups for Some Abelian Examples

The orders of the automorphism groups for $\mathbb{Z}_n \times \mathbb{Z}_2$, $\mathbb{Z}_{p^k}$, and $\mathbb{Z}_p \times \mathbb{Z}_p$ are easily calculable. This section will find these orders, and will consider the general case $Aut(\mathbb{Z}_p \times \mathbb{Z}_p \times \cdots \times \mathbb{Z}_p)$.

### 4.1. The order of $Aut(\mathbb{Z}_n \times \mathbb{Z}_2)$.

Below are calculations supporting the results found in Table 1. The cases for $n$ odd, $n \equiv 0 \pmod{4}$, and $n \equiv 2 \pmod{4}$ will be treated independently.

*Case 1: Let $n$ be odd.* Because $\mathbb{Z}_n \times \mathbb{Z}_2 \cong \mathbb{Z}_{2n}$ (by Theorem 3.3), an application of Theorem 3.7 gives the equality $|Aut(\mathbb{Z}_n \times \mathbb{Z}_2)| = \phi(n)$.

In analyzing $Aut(\mathbb{Z}_n \times \mathbb{Z}_2)$ for $n$ even, let $r$ generate $\mathbb{Z}_n$ and let $s$ generate $\mathbb{Z}_2$ so that $\mathbb{Z}_n \times \mathbb{Z}_2 = \langle r \rangle \times \langle s \rangle$. A homomorphism $\gamma : \mathbb{Z}_n \times \mathbb{Z}_2 \to \mathbb{Z}_n \times \mathbb{Z}_2$ is an element of $Aut(\mathbb{Z}_n \times \mathbb{Z}_2)$ if and only if $\gamma(r)$ has order $n$, $\gamma(s)$ has order 2, and $\gamma(r)^{(n/2)} \neq \gamma(s)$. This is true because a $\gamma(r)$ and $\gamma(s)$ will generate a group $\mathbb{Z}_n \times \mathbb{Z}_2$ only if the intersection of $\langle \gamma(r) \rangle$ and $\langle \gamma(s) \rangle$ is the trivial group $\{e\}$. Therefore finding the order of $Aut(\mathbb{Z}_n \times \mathbb{Z}_2)$ will be a matter of counting the pairs $(x, y)$ that satisfy all of the following:

- $x \in \mathbb{Z}_n \times \mathbb{Z}_2$ has order $n$
- $y \in \mathbb{Z}_n \times \mathbb{Z}_2$ has order $2$
- $x^{(n/2)} \neq y$

*Case 2:* Let $n \equiv 0 \pmod 4$. An element $x$ in $\mathbb{Z}_n \times \mathbb{Z}_2$ has order $n$ if and only if one of the following are satisfied:

- $x = r^i$ and $\gcd(i, n) = 1$
- $x = r^i s$ and $\gcd(i, n) = 1$

Given an element $x \in \mathbb{Z}_n \times \mathbb{Z}_2$ such that $|x| = n$, there are two elements $y$ of order $2$ that satisfy $x^{(n/2)} \neq y$: $s$ and $r^{(n/2)}s$. The element $r^{(n/2)}$ also has order $2$, but $n \equiv 0 \pmod 4$ gives $x^{(n/2)} = r^{(n/2)}$ for every $x$ of order $n$.

Therefore we can find two elements of $\operatorname{Aut}(\mathbb{Z}_n \times \mathbb{Z}_2)$ for each of the $2\phi(n)$ elements $x$ satisfying $|x| = n$: the map sending $r \mapsto x$ and $s \mapsto s$, and the map send $r \mapsto x$ and $s \mapsto r^{(n/2)}s$. Hence we have $|\operatorname{Aut}(\mathbb{Z}_n \times \mathbb{Z}_n)| = 4\phi(n)$.

*Case 3:* Let $n \equiv 2 \pmod 4$. As in the previous case, $r^i$ and $r^i s$ have order $n$ whenever $\gcd(i, n) = 1$. These are not all of the possible generators of $\mathbb{Z}_n$, however: as $\frac{n}{2}$ is an odd number, $|r^j s| = n$ whenever $|r^j| = \frac{n}{2}$. Given that there are $\phi(\frac{n}{2})$ elements $r^j$ that have order $\frac{n}{2}$, a total of $2\phi(n) + \phi(\frac{n}{2})$ elements of order $n$ exist in $\mathbb{Z}_n \times \mathbb{Z}_2$. By Theorem 3.7, $\qquad 2\phi(n) \qquad + \qquad \phi(\frac{n}{2}) \qquad = \qquad 3\phi(n)$.

The elements in $\mathbb{Z}_n \times \mathbb{Z}_2$ with order $2$ are $s$, $r^{(n/2)}$, and $r^{(n/2)}s$. For each elemenet $x$ in $\mathbb{Z}_n \times \mathbb{Z}_2$ such that $|x| = n$, there are two elements $y$ such that $y \neq x^{(n/2)}$.

In conclusion, $|\operatorname{Aut}(\mathbb{Z}_n \times \mathbb{Z}_2)| = 6\phi(n)$ whenever $n \equiv 2 \pmod 4$.

> Several cases for $n$ even:
> $\operatorname{Aut}(\mathbb{Z}_2 \times \mathbb{Z}_2) \cong D_3$,
> $\operatorname{Aut}(\mathbb{Z}_4 \times \mathbb{Z}_2) \cong D_4$,
> $\operatorname{Aut}(\mathbb{Z}_6 \times \mathbb{Z}_2) \cong D_6$,
> $\operatorname{Aut}(\mathbb{Z}_8 \times \mathbb{Z}_2) \cong \mathbb{Z}_2 \times D_4$.

Notice that when $n$ is even, there exists an isomorphism between $\operatorname{Aut}(\mathbb{Z}_n) \times \mathbb{Z}_2$ and the subgroup of $\operatorname{Aut}(\mathbb{Z}_n \times \mathbb{Z}_2)$ consisting of automorphisms that send $s$ to itself. These automorphisms must send $r$ to $r^i$ or $r$ to $r^i s$ for some $i$ relatively prime to $n$.

**Theorem 4.1.** *When $n \in \mathbb{N}$ is a multiple of $4$, $Aut(\mathbb{Z}_n \times \mathbb{Z}_2) \cong \bigl(Aut(\mathbb{Z}_n) \times \mathbb{Z}_2\bigr) \rtimes \mathbb{Z}_2$.*

*Proof.* Let $n \equiv 0 \pmod 4$, and let $[x, y]$ be shorthand notation for the automorphism of $\mathbb{Z}_n \times \mathbb{Z}_2$ that sends $r$ to $x$ and $s$ to $y$ for some $x, y \in \mathbb{Z}_n \times \mathbb{Z}_2$. Given the subgroup $W$ of $\operatorname{Aut}(\mathbb{Z}_n \times \mathbb{Z}_2)$ consisting of automorphisms that send $s$ to itself, we will show first that $\operatorname{Aut}(\mathbb{Z}_n \times \mathbb{Z}_2) \cong W \rtimes \mathbb{Z}_2$ and then that $W \cong \operatorname{Aut}(\mathbb{Z}_n) \times \mathbb{Z}_2$.

The automorphisms in $W$ are of the form $[r^i, s]$ or $[r^i s, s]$ for $i$ relatively prime to $n$. The other elements of $\operatorname{Aut}(\mathbb{Z}_n \times \mathbb{Z}_2)$ are given by $[r^i, r^{(n/2)}s]$ or $[r^i s, r^{(n/2)}s]$ for $\gcd(i, n) = 1$, so by Theorem 3.2 we have $W \triangleleft \operatorname{Aut}(\mathbb{Z}_n \times \mathbb{Z}_2)$. Given the natural injection $\iota : W \to \operatorname{Aut}(\mathbb{Z}_n \times \mathbb{Z}_2)$ as well as the surjection $\tau : \operatorname{Aut}(\mathbb{Z}_n \times \mathbb{Z}_2) \to \mathbb{Z}_2$ that sends elements of $W$ to the identity element in $\mathbb{Z}_2$, we have a short exact sequence $\{e\} \to W \xrightarrow{\iota} \operatorname{Aut}(\mathbb{Z}_n \times \mathbb{Z}_2) \xrightarrow{\tau} \mathbb{Z}_2 \to \{e\}$ that splits via the section sending the nonidentity element in $\mathbb{Z}_2$ to $[r, r^{(n/2)}s]$ in $\operatorname{Aut}(\mathbb{Z}_n \times \mathbb{Z}_2)$. Therefore we have $\operatorname{Aut}(\mathbb{Z}_n \times \mathbb{Z}_2) \cong W \rtimes \mathbb{Z}_2$.

Now we will show that $W \cong \operatorname{Aut}(\mathbb{Z}_n) \times \mathbb{Z}_2$. There is a short exact sequence $\{e\} \to \operatorname{Aut}(\mathbb{Z}_n) \xrightarrow{\alpha} W \xrightarrow{\beta} \mathbb{Z}_2 \to \{e\}$ where $\alpha$ maps $\{r \mapsto r^i\}$ to $[r^i, s]$ and where $\beta$ sends automorphisms of the form $[r^i s, s]$ to the nonidentity element in $\mathbb{Z}_2$. There is a section $\rho : \mathbb{Z}_2 \to W$ given by sending the nonidentity element to $[rs, s]$. Therefore $W \cong \text{Aut}(\mathbb{Z}_n) \rtimes \mathbb{Z}_2$. This group is in fact a direct product because $[r^i, s][rs, s] = [rs, s][r^i, s]$ for every $i$ relatively prime to $n$. $\square$

The theorem above gave an injective map from $\text{Aut}(\mathbb{Z}_n)$ to $\text{Aut}(\mathbb{Z}_n \times \mathbb{Z}_2)$, defined by sending $\{r \mapsto r^i\}$ from $\text{Aut}(\mathbb{Z}_n)$ to $[r^i, s]$ in $\text{Aut}(\mathbb{Z}_n \times \mathbb{Z}_2)$. This function is homomorphic regardless of which natural number $n$ is chosen. Proposition 5.1 in the following section will generalize this notion to inclusion of $\text{Aut}(K)$ as a subgroup of $\text{Aut}(K \times H)$ for any groups $K$ and $H$.

4.2. **The order of $\text{Aut}(\mathbb{Z}_p \times \mathbb{Z}_p \times \cdots \times \mathbb{Z}_p)$.**
Below are explicit calculations for for the orders of $\text{Aut}(\mathbb{Z}_p)$ and $\text{Aut}(\mathbb{Z}_p \times \mathbb{Z}_p)$, as well as a formula for the size of the automorphism group of an arbitrarily long direct product $\mathbb{Z}_p \times \mathbb{Z}_p \times \cdots \times \mathbb{Z}_p$. As discussed in section 2, $p$ represents a prime number and
$$\prod^m \mathbb{Z}_p = \mathbb{Z}_p \times \mathbb{Z}_p \times \cdots \times \mathbb{Z}_p$$
denotes the direct product of $m$ copies of $\mathbb{Z}_p$.

In $\mathbb{Z}_p$, each of the $p - 1$ nonidentity element has order $p$. Letting $r$ generate $\mathbb{Z}_p$, the map $r \mapsto r^i$ is an element of $\text{Aut}(\mathbb{Z}_p)$ whenever $i \in [1, p-1]$ and therefore, $|\text{Aut}(\mathbb{Z}_p)| = \phi(p) = p - 1$.

In calculating the order of $\text{Aut}(\mathbb{Z}_p \times \mathbb{Z}_p)$, let $r$ and $t$ each generate groups of order $p$ so that $\mathbb{Z}_p \times \mathbb{Z}_p = \langle r \rangle \times \langle t \rangle$. A homomorphism $\gamma : \mathbb{Z}_p \times \mathbb{Z}_p \to \mathbb{Z}_p \times \mathbb{Z}_p$ is an automorphism if and only if $|\gamma(r)| = |\gamma(t)| = p$ and $\langle \gamma(r) \rangle$ intersects with $\langle \gamma(s) \rangle$ only at the identity element. To find the order of $\text{Aut}(\mathbb{Z}_p \times \mathbb{Z}_p)$ we must count the pairs $(x, y)$ of elements in $\mathbb{Z}_p \times \mathbb{Z}_p$ such that $\gamma(r) = x$ and $\gamma(t) = y$ determines an automorphism.

Each of the $p^2 - 1$ nonidentity elements in $\mathbb{Z}_p \times \mathbb{Z}_p$ have order $p$, so a given element of $\text{Aut}(\mathbb{Z}_n \times \mathbb{Z}_n)$ may map $r$ to any of $p^2 - 1$ different places. Letting $x$ be a nonidentity element, we must count the elements $y$ of $\mathbb{Z}_p \times \mathbb{Z}_p$ such that $y = p$ and $\langle x \rangle \cap \langle y \rangle = \{e\}$. Each $x$ generates a group of order $p$, and any of the $p^2 - p$ elements of $\mathbb{Z}_p \times \mathbb{Z}_p$ lying outside $\langle x \rangle$ will generate a group of order $p$ that intersects $\langle x \rangle$ only at the identity element. Therefore the order of $\text{Aut}(\mathbb{Z}_p \times \mathbb{Z}_p)$ is $(p^2 - 1)(p^2 - p)$.

Now we will find the order of $\text{Aut}(\mathbb{Z}_p \times \mathbb{Z}_p \times \cdots \times \mathbb{Z}_p)$.
Write $G = \prod^m \mathbb{Z}_p$, and let $\{g_1, g_2, \ldots, g_m\}$ be a set of generators for $G$ so that $\mathbb{Z}_p \times \mathbb{Z}_p \times \cdots \times \mathbb{Z}_p = \langle g_1 \rangle \times \langle g_2 \rangle \times \cdots \times \langle g_m \rangle$. As in $\mathbb{Z}_p \times \mathbb{Z}_p$, every nonidentity element of $G$ has order $p$.

To find the order of $\text{Aut}(G)$, we must count the number of injective maps from the above generators to nonidentity elements that generate groups intersecting only at the identity element. Supposing that an automorphism of $G$ sends $g_1$ to some element $x$ in $G$, there are $p^m - p$ elements $y \in G$ such that $\langle x \rangle \cap \langle y \rangle = \{e\}$. Supposing further that this automorphism is given by $g_1 \mapsto x$ and $g_2 \mapsto y$ for some $y$ not in $\langle x \rangle$, there remain $p^m - p^2$ elements $z \in G$ that are outside of $\langle x \rangle \times \langle y \rangle$. Sending $g_3$ to any such $z$ gives $(\langle x \rangle \times \langle y \rangle) \cap \langle z \rangle = \{e\}$.

Continuing in this fashion, one could specify where an automorphism of $G$ sends the first $n$ generators and then find $p^m - p^n$ elements in $G$ to which the next generator might be sent.

Therefore we can say $\left|\text{Aut}\left(\prod\limits^m \mathbb{Z}_p\right)\right| = \prod\limits_{x=0}^{m-1}(p^m - p^x)$.

**Proposition 4.2.** *Given a prime number $p$ and a natural number $k$, the order of $Aut(\mathbb{Z}_{p^k})$ is $p^k - p^{k-1}$.*

*Proof.* Suppose that $r$ generates $\mathbb{Z}_{p^k}$. There are $\frac{p^k}{p} = p^{k-1}$ elements $r^x$ in $\mathbb{Z}_{p^k}$ such that $p$ divides $x$. This leaves $p^k - p^{k-1}$ elements $r^y$ in $\mathbb{Z}_{p^k}$ such that $\gcd(y, p^k) = 1$. Every one of these elements $r^y$ generates a group of order $p^k$. Therefore $|\text{Aut}(\mathbb{Z}_{p^k})| = p^k - p^{k-1}$. □

## 5. Direct Products

In this section we will give results concerning the automorphism groups for direct products $K \times H$. Although the structure of $\text{Aut}(K \times H)$ may be complicated, we will show that there always exists a subgroup $\text{Aut}(K) \times \text{Aut}(H) < \text{Aut}(K \times H)$. Later we will find sufficient conditions for the isomorphism $\text{Aut}(K) \times \text{Aut}(H) \cong \text{Aut}(K \times H)$.

Notice that $\text{Aut}(K)$ is always isomorphic to a subgroup of $\text{Aut}(K \times H)$: for each $\omega \in \text{Aut}(K)$ there is an element of $\text{Aut}(K \times H)$ that sends $kh$ to $\omega(k)h$ for any given elements $k \in K$ and $h \in H$. These elements of $\text{Aut}(K \times H)$ form a subgroup with the same multiplication as in $\text{Aut}(K)$. This is clearly symmetric, i.e. $\text{Aut}(H) < \text{Aut}(K \times H)$.

**Proposition 5.1.** *$Aut(K \times H)$ has subgroups $Aut(K)$ and $Aut(H)$ that intersect only at the identity automorphism $id_{K \times H}$.*

*Proof.* There exist injective homomorphisms $\zeta : \text{Aut}(K) \to \text{Aut}(K \times H)$ and $\lambda : \text{Aut}(H) \to \text{Aut}(K \times H)$ where $\zeta_\omega$ sends $(k, h)$ to $(\omega(k), h)$ and $\lambda_\delta$ sends $(k, h)$ to $(k, \delta(h))$ for all $\omega \in \text{Aut}(K)$ and $\delta \in \text{Aut}(H)$. The intersection of $\text{im}(\zeta)$ and of $\text{im}(\lambda)$ contains only the automorphism that fixes every element $(k, h)$ of $K \times H$. □

So given a direct product of groups $K \times H$, the automorphism group $\text{Aut}(K \times H)$ has a subgroup $\text{Aut}(K) \times \text{Aut}(H)$ containing precisely those elements of $\text{Aut}(K \times H)$ that send elements from $K$ to $K$ and from $H$ to $H$.

**Corollary 5.2.** *There is an isomorphism $Aut(K \times H) \cong Aut(K) \times Aut(H)$ if and only if both $K$ and $H$ are characteristic in $K \times H$.*

*Proof.* Suppose that $K$ and $H$ are characteristic subgroups of $K \times H$. Then the subgroup $\text{Aut}(K) \times \text{Aut}(H)$ of $\text{Aut}(K \times H)$ given by Proposition 5.1 contains precisely those automorphisms which send elements of $K$ to elements of $K$ and elements of $H$ to elements of $H$. Therefore $\text{Aut}(K) \times \text{Aut}(H)$ must fill out the group $\text{Aut}(K \times H)$.

Now suppose that we are given an isomorphism $\text{Aut}(K \times H) \cong \text{Aut}(K) \times \text{Aut}(H)$. As every element of $\text{Aut}(K \times H)$ will send elements from $K$ to $K$ and elements from $H$ to $H$, the groups $K$ and $H$ must be characteristic in $K \times H$. □

The following two propositions will make direct use of Corollary 5.2.

**Proposition 5.3.** *If $\gcd(m, n) = 1$ then $Aut(\mathbb{Z}_m \times \mathbb{Z}_n) \cong Aut(\mathbb{Z}_m) \times Aut(\mathbb{Z}_n)$.*

*Proof.* Suppose that $\gcd(m,n) = 1$. We can apply Corollary 5.2 after showing that $\mathbb{Z}_m$ and $\mathbb{Z}_n$ are characteristic subgroups of $\mathbb{Z}_m \times \mathbb{Z}_n$. Writing $\mathbb{Z}_m \times \mathbb{Z}_n = \langle a \rangle \times \langle b \rangle$, each automorphism of $\mathbb{Z}_m \times \mathbb{Z}_n$ can be given by where it sends $a$ and $b$. Any such automorphism must send $a$ to an element of order $m$. Given that $|a^i b^j| = \text{lcm}(|a^i|, |b^j|)$, we have

$$\text{lcm}(|a^i|, |b^j|) = m \iff |a^i| = m \text{ and } |b^j| = 1$$

because the order of $b$ is relatively prime to $m$. Therefore, every automorphism of $\mathbb{Z}_m \times \mathbb{Z}_n$ must send $a$ to another generator of $\langle a \rangle$. A similar argument gives that each element of $\text{Aut}(\mathbb{Z}_m \times \mathbb{Z}_n)$ must send $b$ to a power of $\langle b \rangle$. In conclusion, the automorphism of $\mathbb{Z}_m \times \mathbb{Z}_n$ are of the form $\begin{pmatrix} a \mapsto \chi(a) \\ b \mapsto \delta(b) \end{pmatrix}$ for $\chi \in \text{Aut}(\mathbb{Z}_m)$ and $\delta \in \text{Aut}(\mathbb{Z}_n)$. $\square$

Theorem 3.3 lets us decompose any group $\mathbb{Z}_n$ into a direct product of cyclic groups with relatively prime order, and Proposition 5.3 allows us to write

$$\text{Aut}(\mathbb{Z}_n) \cong \text{Aut}(\mathbb{Z}_{p_1^{k_1}}) \times \text{Aut}(\mathbb{Z}_{p_2^{k_2}}) \times \cdots \times \text{Aut}(\mathbb{Z}_{p_1^{k_1}})$$

assuming that $n = p_1^{k_1} p_2^{k_2} \cdots p_w^{k_w}$ for $p_1, p_2, \ldots, p_w$ distinct primes. Thus the result from Proposition 4.2 gives

$$|\text{Aut}(\mathbb{Z}_n)| = (p_1^{k_1} - p_1^{k_1-1})(p_2^{k_2} - p_2^{k_2-1}) \cdots (p_w^{k_w} - p_w^{k_w-1})$$

for the order of $\text{Aut}(\mathbb{Z}_n)$.

We reached conclusion in Proposition 5.3 because $|a^i|$ and $|b^j|$ divide the order of an element $a^i b^j$ in $\mathbb{Z}_m \times \mathbb{Z}_n \cong \langle a, b \mid a^m = b^n = e \rangle$. The following proposition will generalize this notion to groups of the form $K \times H$.

**Proposition 5.4.** *Let $K$ and $H$ be groups. If $\gcd(|K|, |H|) = 1$ then $Aut(K \times H) \cong Aut(K) \times Aut(H)$.*

*Proof.* The orders of $K$ and of $H$ are relatively prime if and only if $\gcd(|k|, |h|) = 1$ for all $k \in K$ and $h \in H$. This can be shown by using Cauchy's Theorem as well as the fact that the order of each element of a given group divides the order of that group.

Suppose that $\gcd(|k|, |h|) = 1$ for all elements $k \in K$ and $h \in H$. We will show that both $K$ and $H$ are characteristic in $K \times H$ by considering the orders of the elements in $K \times H$.

An automorphism of $K \times H$ must send each $k \in K$ to an element of order $|k|$, but every element $h \in H$ has order relatively prime to $|k|$. Thus given any $x \in K$, each nonidentity element $h \in H$ gives an element $xh \in K \times H$ with order not equal to $|k|$. Similarly, given any $y \in H$, every nonidentity element $k \in K$ gives an element $ky$ with order not equal to $|h|$.

Therefore, given $k \in K$ and $h \in H$, every automorphism of $K \times H$ must be of the form $kh \mapsto \chi(k)\delta(h)$ for some $\chi \in \text{Aut}(K)$ and $\delta \in \text{Aut}(H)$. This gives $\text{Aut}(K \times H) \cong \text{Aut}(K) \times \text{Aut}(H)$. $\square$

## 6. Semidirect Products

It is not as easy to find homomorphisms concerning $K \rtimes_\psi H$ as it is for a direct product $K \times H$. For example, it is possible that given some $\omega \in \mathrm{Aut}(K)$ we might be unable to create an automorphism for $K \rtimes_\psi H$ by sending $kh$ to $\omega(k)h$ for every $k \in K$ and $h \in H$. The following lemma will give a clear way to determine whether such functions as $(kh \mapsto \omega(k)h)$ are homomorphic.

**Lemma 6.1.** *Let $K, H, L$ and $M$ be groups. Let $\gamma \in \hom(K, L), \varphi \in \hom(H, M)$, and $\psi \in \hom\big(H, Aut(K)\big)$. If $F : K \rtimes_\psi H \to L \rtimes M$ is a function defined by $kh \mapsto \gamma(k)\varphi(h)$ for each element $k \in K$ and $h \in H$, then $F$ is homomorphic if and only if $F(hk) = F(h)F(k)$ for all $k \in K$ and $h \in H$.*

*Proof.* If $F(hk) \neq F(h)F(k)$ for some $k \in K$ and $h \in H$, then $F$ is by definition not homomorphic. Now suppose that $F(hk) = F(h)F(k)$ for all $k \in K$ and for all $h \in H$. We want to show that $F(xy) = F(x)F(y)$ for all $x, y \in K \rtimes_\psi H$. Without loss of generality, write $x = kh$ and $y = k'h'$ for some $k, k' \in K$ and $h, h' \in H$. Now we check that $F$ is homomorphic:

$$\begin{aligned}
F(khk'h') = & \ F\big(k\psi_h(k')hh'\big) & =\gamma\big(k\psi_h(k')\big)\varphi(hh') & =\gamma(k)\gamma\big(\psi_h(k')\big)\varphi(h)\varphi(h') \\
= & \ \gamma(k)F\big(\psi_h(k')h\big)\varphi(h') & =\gamma(k)F(hk')\varphi(h') & = \gamma(k)F(h)F(k')\varphi(h') \\
= & \ \gamma(k)\varphi(h)\gamma(k')\varphi(h') & = F(kh)F(k'h')
\end{aligned}$$

$\square$

### 6.1. Generalizations to the semidirect product case.

Having stated Lemma 6.1, we can proceed towards finding sufficient conditions for $\mathrm{Aut}(K) < \mathrm{Aut}(K \rtimes_\psi H)$ and $\mathrm{Aut}(H) < \mathrm{Aut}(K \rtimes_\psi H)$. Let us begin by applying the principle used in Proposition 5.1, defining for each $\omega \in \mathrm{Aut}(K)$ a map $\zeta_\omega$ such that $\zeta_\omega(kh) = \omega(k)h$ for all $k \in K$ and $h \in H$. By Lemma 6.1 we have that $\zeta_\omega$ is homomorphic if and only if $\zeta_\omega(hk) = \zeta_\omega(h)\zeta_\omega(k)$ for all $h \in H$ and $k \in K$. The following theorem gives a condition necessary and sufficient to include $\zeta_\omega$ in $\mathrm{Aut}(K \rtimes_\psi H)$ for all $\omega \in \mathrm{Aut}(K)$.

**Theorem 6.2.** *Let $\psi \in \hom\big(H, Aut(K)\big)$. If $im(\psi) < Z\big(Aut(K)\big)$ then $Aut(K) < Aut(K \rtimes_\psi H)$.*

*Proof.* Suppose that $\mathrm{im}(\psi) < Z\big(\mathrm{Aut}(K)\big)$, and let $\zeta : \mathrm{Aut}(K) \to \mathrm{Aut}(K \rtimes_\psi H)$ send each $\omega \in \mathrm{Aut}(K)$ to the map $\{kh \mapsto \omega(k)h\}$ where $k \in K$ and $h \in H$. Letting $\zeta_\omega$ denote this image of $\omega$ through $\zeta$,

$$\zeta_\omega(hk) = \zeta_\omega(\psi_h(k)h) = \omega(\psi_h(k))h = \psi_h(\omega(k))h = h\omega(k) = \zeta_\omega(h)\zeta_\omega(k)$$

gives that each $\zeta_\omega$ is a homomorphism. Moreover, $\zeta$ is itself a homomorphism because $\zeta_{(\kappa\upsilon)} = \zeta_\kappa \zeta_\upsilon$ for all $\kappa, \upsilon \in \mathrm{Aut}(K)$. Therefore $\mathrm{Aut}(K) \cong \mathrm{im}(\zeta) < \mathrm{Aut}(K \rtimes_\psi H)$. $\square$

Having given a necessary and sufficient condition for application of the map $\zeta$ (from Proposition 5.1) to the semidirect product case, we will consider the map $\lambda$ (from the same theorem) in the semidirect product case as well. Define for each $\delta \in \mathrm{Aut}(H)$ a function $\lambda_\delta$ such that $\lambda_\delta(kh) = k\delta(h)$ for all $k \in K$ and $h \in H$. To demonstrate homomorphism of $\lambda_\delta$, Lemma 6.1 requires that $\lambda_\delta(hk) = \lambda_\delta(h)\lambda_\delta(k)$ for every $h \in H$ and $k \in K$. The following theorem's hypothesis is necessary and sufficient for every $\delta \in \mathrm{Aut}(H)$ to provide an element $\lambda_\delta$ in $\mathrm{Aut}(K \rtimes_\psi H)$.

**Theorem 6.3.** *Let $\psi \in \hom(H, Aut(K))$. If $\psi \circ \delta = \psi$ for every $\delta \in Aut(H)$ then $Aut(H) < Aut(K \rtimes_\psi H)$.*

*Proof.* Let $\psi \circ \delta = \psi$ for all $\delta \in Aut(H)$, and let $\lambda : Aut(H) \to Aut(K \rtimes_\psi H)$ send each $\delta \in Aut(H)$ to the map $\{kh \mapsto k\delta(h)\}$ where $k \in K$ and $h \in H$. Letting $\lambda_\delta$ denote this image of $\delta$ in through $\lambda$,

$$\lambda_\delta(hk) = \lambda_\delta(\psi_h(k)h) = \psi_h(k)\delta(h) = (\psi \circ \delta)_h(k)\delta(h)$$
$$= \psi_{\delta(h)}(k)\delta(h) = \delta(h)k = \lambda_\delta(h)\lambda_\delta(k)$$

gives that each $\lambda_\delta$ is a homomorphism. Moreover, $\lambda$ is itself a homomorphism because $\lambda_{\eta\mu} = \lambda_\eta \lambda_\mu$ for all $\eta, \mu \in Aut(K)$. Therefore $Aut(H) \cong \text{im}(\lambda) < Aut(K \rtimes_\psi H)$. □

In Proposition 5.1 we find subgroups $Aut(K)$ and $Aut(H)$ of $Aut(K \times H)$ given by $\text{im}(\zeta)$ and $\text{im}(\lambda)$, and we see that the intersection of these subgroups contains only $\text{id}_{K \times H}$. We can complete our extension of Proposition 5.1 to the semidirect product case by noting that if $\psi$ satisfies the hypotheses for Theorems 6.2 and 6.3, then the subgroups $Aut(K)$ and $Aut(H)$ of $Aut(K \rtimes_\psi H)$ given by given by $\text{im}(\zeta)$ and $\text{im}(\lambda)$ intersect only at $\text{id}_{K \rtimes_\psi H}$. This is because if $\zeta_\omega = \lambda_\delta$ then $\omega(k)h = k\lambda(h)$ for all $k$ in $K$ and $h$ in $H$, which can occur only if $\omega = \text{id}_K$ and $\delta = \text{id}_H$.

This allows us to obtain an alternate proof for Proposition 5.1 via application of Theorems 6.2 and 6.3: If $K \times H = K \rtimes_\psi H$ then $\psi$ must send every element of $H$ to $\text{id}_K$, giving both $\text{im}(\psi) < Z(Aut(K))$ and $\psi \circ \delta = \psi$ for all $\delta \in Aut(H)$.

Proposition 5.3 does not generalize to the semidirect product case as smoothly as does Proposition 5.1. A weaker verson of Proposition 5.3 is given below.

**Theorem 6.4.** *If $\gcd(m, n) = 1$ then $\mathbb{Z}_m$ is characteristic in $\mathbb{Z}_m \rtimes_\psi \mathbb{Z}_n$ for every $\psi \in \hom(\mathbb{Z}_n, Aut(\mathbb{Z}_m))$.*

*Proof.* Suppose that $\gcd(m, n) = 1$ and let

$$\mathbb{Z}_m \rtimes_\psi \mathbb{Z}_n \cong \langle a, b \mid a^m = b^n = e, \ ba = \psi_b(a)b \rangle$$

for some $\psi \in \hom(\mathbb{Z}_n, Aut(\mathbb{Z}_m))$. We will show that $\gamma(a) = a^i b^j$ implies $j = 0$ for every automorphism $\gamma \in Aut(\mathbb{Z}_m \rtimes_\psi \mathbb{Z}_n)$.

Since $|a| = m$, we know that the order of $a^i b^j$ must equal $m$ as well. Writing $(a^i b^j)^m$ explicitly as

$$(a^i b^j)^m = \left(a^i \psi_b^j(a^i) \psi_b^{2j}(a^i) \cdots \psi_b^{(m-1)j}(a^i)\right) b^{mj}$$

makes clear that if $(a^i b^j)^m = e$ then $b^{mj} = e$. Therefore, if $a$ maps through an automorphism to $a^i b^j$ then the order of $b^{mj}$ must be a multiple of $n$. But $|b^{mj}|$ divides $n$ and $\gcd(m, n) = 1$, so we can conclude that $j = 0$. □

The previous theorem gives that if $\gcd(m, n) = 1$, then both $\mathbb{Z}_m$ and $\mathbb{Z}_n$ are characteristic in $\mathbb{Z}_m \times \mathbb{Z}_n$. Applying the result from Corollary 5.2 gives another proof for Proposition 5.3.

**Theorem 6.5.** *Let $K$ and $H$ be groups. If $\gcd(|K|, |H|) = 1$ then $K$ is characteristic in $K \rtimes_\psi H$ for all $\psi \in \hom(H, Aut(K))$.*

*Proof.* As mentioned in Proposition 5.4, the orders of $K$ and of $H$ are relatively prime if and only if $\gcd(|k|, |h|) = 1$ for all $k \in K$ and $h \in H$. Suppose that indeed

$\gcd(|k|, |h|) = 1$ for all elements $k \in K$ and $h \in H$. We can write
$$(kh)^m = \left(k\psi_h(k)\psi_h^2(k) \cdots \psi_h^{(m-1)}(k)\right) h^m$$
which demonstrates that $|kh|$ is a multiple of $|h|$ for all $k \in K$ and $h \in H$. If $h \in H$ is not the identity element, then the order of $kh$ does not divide $|K|$. Thus every automorphism of $K \rtimes_\psi H$ must send the elements of $K$ amongst themselves. $\square$

## 6.2. Isomorphism of different semidirect products $K \rtimes H$.

This subsection will find a sufficient condition for isomorphism of $K \rtimes_\psi H$ and $K \rtimes_\varphi H$ where $\psi$ and $\varphi$ are elements of $\hom(H, \text{Aut}(K))$.

Let $\mathscr{R}$ be the equivalence relation on $\hom(H, \text{Aut}(K))$ defined by $\psi \sim \varphi$ if and only if $\psi = \varphi \circ \delta$ for some $\delta \in \text{Aut}(H)$. Thus the equivalence class of each $x \in \hom(H, \text{Aut}(K))$ is $\{x \circ y \mid y \in \text{Aut}(H)\}$.

**Theorem 6.6.** *Let $\psi, \varphi \in \hom(H, Aut(K))$. If $\psi \sim \varphi$ then $K \rtimes_\psi H \cong K \rtimes_\varphi H$.*

*Proof.* In $K \rtimes_\psi H$ each element can be expressed as the product $kh$ for some $k \in K$ and $h \in H$. Because $kh$ is also an element of $K \rtimes_\varphi H$, care must be taken not to confuse the elements of $K \rtimes_\psi H$ and of $K \rtimes_\varphi H$.

Suppose that $\psi \sim \varphi$ and that $\delta \in \text{Aut}(H)$ satisfies $\psi = \varphi \circ \delta$. Let the function $\Theta : K \rtimes_\psi H \to K \rtimes_\varphi H$ be defined by $\Theta(kh) \mapsto k\delta(h)$. This function is surjective because for each $kh \in K \rtimes_\varphi H$ we have $\Theta(k\delta^{-1}(h)) = kh$. As the orders of $K \rtimes_\psi H$ and of $K \rtimes_\varphi H$ are equal, we have bijection.

We can check that $\Theta$ is isomorphic by applying Lemma 6.1. Observing that $\psi_h = (\varphi \circ \delta)_h = \varphi_{\delta(h)}$ for all $h \in H$, we find that
$$\Theta(hk) = \Theta(\psi_h(k)h) = \psi_h(k)\delta(h) = \varphi_{\delta(h)}(k)\delta(h) = \delta(h)k = \Theta(h)\Theta(k)$$
for all $k \in K$ and $h \in H$. $\square$

**Corollary 6.7.** *If $H \cong Aut(K)$ and if $\psi, \varphi \in \hom(H, Aut(K))$ are isomorphisms then $K \rtimes_\psi H \cong K \rtimes_\varphi H$.*

*Proof.* Suppose that $H \cong \text{Aut}(K)$ and that $\psi, \varphi \in \hom(H, \text{Aut}(K))$ are isomorphisms. By Theorem 6.6, it will be sufficient to show that $\psi \sim \varphi$. Indeed, we have $\psi = \varphi \circ \delta$ for $\varphi^{-1} \circ \psi = \delta \in \text{Aut}(H)$. $\square$

## 7. Application to Dihedral Groups

The automorphism group $\text{Aut}(G)$ has structure dependant on that of the given group $G$, yet it may prove difficult to give an explicit representation of $\text{Aut}(G)$ even when everything is know about $G$. Finding such a representation is helpful in deciding whether there exists an isomorphism between $G$ and $\text{Aut}(G)$. In this section we will uncover in full detail the structure of the $\text{Aut}(D_n)$ groups, concluding with the result $D_n \cong \text{Aut}(D_n) \iff \phi(n) = 2$.

A healthy first step towards describing $\text{Aut}(D_n)$ is the observation that $\text{Aut}(\mathbb{Z}_n) < \text{Aut}(D_n)$, which can be seen via an application of Theorem 6.2. The following lemma will allow us to apply that Theorem.

**Lemma 7.1.** *For any abelian group $A$, the element of $Aut(A)$ sending each $a \in A$ to $a^{-1}$ is a member of $Z(Aut(A))$.*

*Proof.* Let $\theta \in \mathrm{Aut}(A)$ be defined by $\theta(a) = a^{-1}$ for all $a \in A$. Given any $\gamma \in \mathrm{Aut}(A)$ we have

$$(\theta\gamma\theta)(a) = (\theta\gamma)(a^{-1}) = \theta\big(\gamma(a^{-1})\big) = \theta\Big(\big(\gamma(a)\big)^{-1}\Big) = \gamma(a)$$

for all $a \in A$. □

Letting $\theta \in \mathrm{Aut}(\mathbb{Z}_n)$ be the automorphism that sends each element to its inverse (as above), we can write $D_n \cong \mathbb{Z}_n \rtimes_\kappa \mathbb{Z}_2$ where $\kappa \in \hom\big(\mathbb{Z}_2, \mathrm{Aut}(\mathbb{Z}_n)\big)$ sends the nonidentity element in $\mathbb{Z}_2$ to $\theta \in \mathrm{Aut}(\mathbb{Z}_n)$. Furthermore, let $r$ be a generator for $\mathbb{Z}_n$, and let $s$ generate $\mathbb{Z}_2$ so that $D_n \cong \langle r \rangle \rtimes_\kappa \langle s \rangle$.

The following Theorem describes $\mathrm{Aut}(D_n)$ fully.

**Theorem 7.2.** *Let $n \geq 3$. Then $Aut(D_n) \cong \mathbb{Z}_n \rtimes_\psi Aut(\mathbb{Z}_n)$ where $\psi : Aut(\mathbb{Z}_n) \to Aut(\mathbb{Z}_n)$ is the identity map.*

*Proof.* Using the notation $D_n \cong \langle r \rangle \rtimes_\kappa \langle s \rangle$ we will show that $\mathbb{Z}_n \triangleleft \mathrm{Aut}(D_n)$ and that the short exact sequence $\{e\} \to \mathbb{Z}_n \to \mathrm{Aut}(D_n) \to \mathrm{Aut}(\mathbb{Z}_n) \to \{e\}$ splits. Later we will demonstrate that if a homomorphism $\psi$ satisfies $\mathrm{Aut}(D_n) \cong \mathbb{Z}_n \rtimes_\psi \mathrm{Aut}(\mathbb{Z}_n)$, then $\psi$ must be a bijection.

In the context of this proof, let $[x, y]$ denote the automorphism of $D_n$ given by $r \mapsto x$ and $s \mapsto y$. The element, $[r, rs] \in \mathrm{Aut}(D_n)$ represents the map sending each $r^i$ to itself and each $r^i s$ to $r^{(i+1)}s$. This map is an automorphism because $\langle rs \rangle \cong \mathbb{Z}_2$ and conjugation of any $r^i$ by $rs$ yields $r^{-i}$. It so happens that $[r, rs]^2 = [r, r^2 s]$, and more generally $[r, r^u s][r, r^v s] = [r, r^{(u+v)} s]$. Therefore the cyclic subgroup of $\mathrm{Aut}(D_n)$ generated by $[r, rs]$ has order $n$, and the elements of this subgroup take the form $[r, r^j s]$ for integers $j$.

We will now find the general form for automorphisms of $D_n$: $r$ must always map to an element of order $n$, and all of these elements live in $\langle r \rangle$. By Lemma 7.1, $\kappa(\mathbb{Z}_2) < Z\big(\mathrm{Aut}(\mathbb{Z}_n)\big)$, and therefore Theorem 6.2 supplies a unique element $\zeta_\omega \in \mathrm{Aut}(D_n)$ for each $\omega \in \mathrm{Aut}(\mathbb{Z}_n)$ given by $\zeta_\omega = [\omega(r), s]$. Composing the automorphisms $[\omega(r), s]$ and $[r, r^j s]$ gives $[r, r^j s][\omega(r), s] = [\omega(r), r^j s]$. Because $s$ is not contained in $\langle r \rangle$, every automorphism $\gamma$ of $D_n$ must send $s$ to an an element of order two that is not contained in $\langle \gamma(r) \rangle$. The subgroup $\langle r \rangle$ of $D_n$ is characteristic because $\gamma(r)$ must have order $n$, so $\gamma$ send $s$ to some $r^j s \in D_n$. Thus we know that every automorphism of $D_n$ must be of the form $[\omega(r), r^j s]$ for some $\omega \in \mathrm{Aut}(Z_n)$.

The subgroup $\langle [r, rs] \rangle \cong \mathbb{Z}_n$ of $\mathrm{Aut}(D_n)$ is normal: conjugation of $[r, rs]$ by $[\omega(r), s]$ gives $[r, \omega(r)s]$. The quotient group $\mathrm{Aut}(D_n)/\langle [r, rs] \rangle$ is isomorphic to $\mathrm{Aut}(\mathbb{Z}_n)$, and a section that splits the sequence

$$\{e\} \to \mathbb{Z}_n \to \mathrm{Aut}(D_n) \to \mathrm{Aut}(\mathbb{Z}_n) \to \{e\}$$

can be given by $\zeta : \mathrm{Aut}(\mathbb{Z}_n) \to \mathrm{Aut}(D_n)$ which sends each $\omega \in \mathrm{Aut}(\mathbb{Z}_n)$ to $[\omega(r), s]$. Therefore $\mathrm{Aut}(D_n) \cong \mathbb{Z}_n \rtimes_\psi \mathrm{Aut}(\mathbb{Z}_n)$ where $\psi \in \hom\big(\mathrm{Aut}(\mathbb{Z}_n), \mathrm{Aut}(\mathbb{Z}_n)\big)$.

This homomorphism $\psi$ is in fact the identity map. We know this because conjugation of $[r, rs]$ by $[\omega(r), s]$ gives $[r, \omega(r)s]$. □

This description of the $\mathrm{Aut}(D_n)$ groups does not leave much to the imagination. Note that by Corollary 6.7, every isomorphism $\psi \in \hom\big(\mathrm{Aut}(\mathbb{Z}_n), \mathrm{Aut}(\mathbb{Z}_n)\big)$ supplies the same group $\mathbb{Z}_n \rtimes_\psi \mathrm{Aut}(\mathbb{Z}_n) \cong D_n$. Now it is possible to deduce precisely when $D_n \cong \mathrm{Aut}(D_n)$.

**Corollary 7.3.** $Aut(D_n) \cong D_n$ *if and only if* $\phi(n) = 2$.

*Proof.* When $n = 1$ or $n = 2$ we have $\phi(n) = 1$, which gives $D_1 \cong \mathbb{Z}_2 \not\cong \{e\} \cong \mathrm{Aut}(D_1)$ and $D_2 \cong \mathbb{Z}_2 \times \mathbb{Z}_2 \not\cong D_3 \cong \mathrm{Aut}(\mathbb{Z}_2 \times \mathbb{Z}_2)$.

Now let $n \geq 3$. We will show that given an isomorphism $\psi$ in $\hom\bigl(\mathrm{Aut}(\mathbb{Z}_n), \mathrm{Aut}(\mathbb{Z}_n)\bigr)$, we can say $\mathbb{Z}_n \rtimes_\kappa \mathbb{Z}_2 \cong \mathbb{Z}_n \rtimes_\psi \mathrm{Aut}(Z_n)$ if and only if $\mathbb{Z}_2 \cong \mathrm{Aut}(\mathbb{Z}_n)$.

Suppose that $\phi(n) = 2$, and notice that if $D_n \cong \mathbb{Z}_n \rtimes_\kappa \mathbb{Z}_2$ then $\kappa : \mathbb{Z}_2 \to \mathrm{Aut}(\mathbb{Z}_n)$ is an isomorphism. Invoking Corollary 6.7 we have $\mathbb{Z}_n \rtimes_\kappa \mathbb{Z}_2 \cong \mathbb{Z}_n \rtimes_\psi \mathrm{Aut}(\mathbb{Z}_n)$. To find that $D_n \cong \mathrm{Aut}(D_n)$ implies $\phi(n) = 2$, observe that the groups $D_n$ and $\mathrm{Aut}(D_n)$ have orders $2n$ and $n\phi(n)$, respecively. These orders are equal only if $\phi(n) = 2$. □

## 8. $\mathrm{Aut}(\mathbb{Z}_8 \rtimes \mathbb{Z}_2)$ for each action of $\mathbb{Z}_2$ on $\mathbb{Z}_8$

In Section 4 the $\mathrm{Aut}(\mathbb{Z}_n \times \mathbb{Z}_2)$ groups were classified by order, and in Section 7 the structures of the $\mathrm{Aut}(D_n)$ groups were given. This section will focus on the four $\mathrm{Aut}(\mathbb{Z}_8 \rtimes \mathbb{Z}_2)$ groups, providing a glimpse of what lies beyond the abelian and dihedral cases for $\mathrm{Aut}(\mathbb{Z}_m \rtimes \mathbb{Z}_n)$. It will be shown that $\mathrm{Aut}(\mathbb{Z}_8 \rtimes \mathbb{Z}_2) \cong \mathbb{Z}_2 \times D_4$ whenever $\mathbb{Z}_8 \rtimes \mathbb{Z}_2$ is not dihedral, and that $\mathrm{Aut}(D_8) \cong (\mathbb{Z}_2 \times D_4) \rtimes \mathbb{Z}_2$.

The cyclic group $\mathbb{Z}_8$ is the smallest cyclic group that has more than one automorphism of order 2; for $n \leq 7$ we have either $\mathbb{Z}_n \rtimes \mathbb{Z}_2 \cong \mathbb{Z}_n \times \mathbb{Z}_2$ or $\mathbb{Z}_n \rtimes \mathbb{Z}_2 \cong D_n$. As $\mathrm{Aut}(\mathbb{Z}_8) \cong \mathbb{Z}_2 \times \mathbb{Z}_2$, each automorphism of $\mathbb{Z}_8$ produces a unique element of $\hom\bigl(\mathbb{Z}_2, \mathrm{Aut}(\mathbb{Z}_8)\bigr)$ and a unique group $\mathbb{Z}_8 \rtimes \mathbb{Z}_2$.

8.1. **The actions of $\mathbb{Z}_2$ on $\mathbb{Z}_8$.** Let $r$ and $s$ generate $\mathbb{Z}_8$ and $\mathbb{Z}_2$, respectively. There are four automorphisms in $\mathrm{Aut}(\mathbb{Z}_8)$, expressed here by where they send $r$:

$$\begin{array}{ll} (r \mapsto r) & (r \mapsto r^3) \\ (r \mapsto r^5) & (r \mapsto r^7) \end{array}$$

Notice that the three nonidentity elements of this automorphism group have order 2, and that the product of any two nonidentity elements is equal to the third nonidentity element. We may therefore write that $\mathrm{Aut}(\mathbb{Z}_8) \cong \mathbb{Z}_2 \times \mathbb{Z}_2$.

The set $\hom\bigl(\mathbb{Z}_2, \mathrm{Aut}(\mathbb{Z}_8)\bigr)$ has 4 elements–one for each automorphism in $\mathrm{Aut}(\mathbb{Z}_8)$. These homomorphisms from $\mathbb{Z}_2$ to $\mathrm{Aut}(\mathbb{Z}_8)$ can be given by

$$\begin{aligned} \rho(s) &= (r \mapsto r) \\ \sigma(s) &= (r \mapsto r^3) \\ \tau(s) &= (r \mapsto r^5) \\ \upsilon(s) &= (r \mapsto r^7). \end{aligned}$$

*Remark* 8.1. $\mathbb{Z}_8 \rtimes_\rho \mathbb{Z}_2 \cong \mathbb{Z}_8 \times \mathbb{Z}_2$ and $\mathbb{Z}_8 \rtimes_\upsilon \mathbb{Z}_2 \cong D_8$, while $sr = r^3 s$ in $\mathbb{Z}_8 \rtimes_\sigma \mathbb{Z}_2$ and $sr = r^5 s$ in $\mathbb{Z}_8 \rtimes_\tau \mathbb{Z}_2$

8.2. **General forms of the automorphisms.** For a given group $G = \mathbb{Z}_8 \rtimes \mathbb{Z}_2$, let $[r^a s^b, r^c s^d]$ denote the automorphism of $G$ that sends $r^i s^j \mapsto (r^a s^b)^i (r^c s^d)^j$ for all $r^i s^j \in G$. Multiplication of elements in $\mathrm{Aut}(G)$ is given by $[r^a s^b, r^c s^d][r^w s^x, r^y s^z] = [(r^a s^b)^w (r^c s^d)^x, (r^a s^b)^y (r^c s^d)^z]$. The following lines describe in general the elements of the automorphism groups considered.

$$\begin{aligned} x \in \mathrm{Aut}(\mathbb{Z}_8 \times \mathbb{Z}_2) &\iff x = [r^i s^j, r^{4k} s] \text{ for some } i \in \{1,3,5,7\} \text{ and } j,k \in \{0,1\}. \\ x \in \mathrm{Aut}(\mathbb{Z}_8 \rtimes_\sigma \mathbb{Z}_2) &\iff x = [r^i, r^k s] \text{ for some } i \in \{1,3,5,7\} \text{ and } k \in \{0,2,4,6\}. \\ x \in \mathrm{Aut}(\mathbb{Z}_8 \rtimes_\tau \mathbb{Z}_2) &\iff x = [r^i s^j, r^{4k} s] \text{ for some } i \in \{1,3,5,7\} \text{ and } j,k \in \{0,1\}. \\ x \in \mathrm{Aut}(D_8) &\iff x = [r^i, r^k s] \text{ for some } i \in \{1,3,5,7\} \text{ and } k \in \{0,1,2,3,4,5,6,7\}. \end{aligned}$$

8.3. **The structure of $\text{Aut}(\mathbb{Z}_8 \rtimes \mathbb{Z}_2)$ for $\mathbb{Z}_8 \rtimes \mathbb{Z}_2 \not\cong D_8$.** Even though the elements of $\text{Aut}(\mathbb{Z}_8 \rtimes \mathbb{Z}_2)$ look different for each action of $\mathbb{Z}_2$ on $\mathbb{Z}_8$, very similar arguments can be used to expose the structures of $\text{Aut}(\mathbb{Z}_8 \rtimes_\rho \mathbb{Z}_2)$, $\text{Aut}(\mathbb{Z}_8 \rtimes_\sigma \mathbb{Z}_2)$, and $\text{Aut}(\mathbb{Z}_8 \rtimes_\tau \mathbb{Z}_2)$.

**Theorem 8.2.** $Aut(\mathbb{Z}_8 \times \mathbb{Z}_2) \cong \mathbb{Z}_2 \times D_4$.

*Proof.* The element $[r^3, s]$ is central in $\text{Aut}(\mathbb{Z}_8 \times \mathbb{Z}_2)$, and therefore $\langle [r^3, s] \rangle \cong \mathbb{Z}_2$ is normal in $\text{Aut}(\mathbb{Z}_8 \times \mathbb{Z}_2)$. Letting $Q$ denote the quotient $\text{Aut}(\mathbb{Z}_8 \times \mathbb{Z}_2)/\langle [r^3, s] \rangle$, the short exact sequence

$$\{e\} \longrightarrow \langle [r^3, s] \rangle \longrightarrow \text{Aut}(\mathbb{Z}_8 \times \mathbb{Z}_2) \longrightarrow Q \longrightarrow \{e\}$$

splits, giving $\text{Aut}(\mathbb{Z}_8 \times \mathbb{Z}_2) \cong \mathbb{Z}_2 \times Q$.

The image of $Q$ in $\text{Aut}(\mathbb{Z}_8 \times \mathbb{Z}_2)$ through one possible section gives a subgroup containing elements of the form $[r^i s^j, r^{4k} s]$ for $i \in \{1, 5\}$ and $j, k \in \{0, 1\}$. The group of order four generated by $[rs, r^4 s]$ is normal in this subgroup, and the short exact sequence

$$\{e\} \longrightarrow \langle [rs, r^4 s] \rangle \longrightarrow Q \longrightarrow \langle [r, r^4 s] \rangle \longrightarrow \{e\}$$

splits, demonstrating that $Q$ is isomorphic to $\mathbb{Z}_4 \rtimes \mathbb{Z}_2$. Conjugation of $[rs, r^4 s]$ by $[r, r^4 s]$ yields $[rs, r^4 s]^{-1}$, and therefore $Q \cong D_4$. In conclusion, $\text{Aut}(\mathbb{Z}_8 \times \mathbb{Z}_2) \cong \mathbb{Z}_2 \times D_4$. $\square$

**Theorem 8.3.** $Aut(\mathbb{Z}_8 \rtimes_\sigma \mathbb{Z}_2) \cong \mathbb{Z}_2 \times D_4$.

*Proof.* The element $[r^5, s]$ is central in $\text{Aut}(\mathbb{Z}_8 \rtimes_\sigma \mathbb{Z}_2)$, and therefore $\langle [r^5, s] \rangle \cong \mathbb{Z}_2$ is normal in $\text{Aut}(\mathbb{Z}_8 \rtimes_\sigma \mathbb{Z}_2)$. Letting $U$ denote the quotient $\text{Aut}(\mathbb{Z}_8 \rtimes_\sigma \mathbb{Z}_2)/\langle [r^5, s] \rangle$, the short exact sequence

$$\{e\} \longrightarrow \langle [r^5, s] \rangle \longrightarrow \text{Aut}(\mathbb{Z}_8 \rtimes_\sigma \mathbb{Z}_2) \longrightarrow U \longrightarrow \{e\}$$

splits, giving $\text{Aut}(\mathbb{Z}_8 \rtimes_\sigma \mathbb{Z}_2) \cong \mathbb{Z}_2 \times U$.

The image of $U$ in $\text{Aut}(\mathbb{Z}_8 \rtimes_\sigma \mathbb{Z}_2)$ through one possible section gives a subgroup containing elements of the form $[r^i, r^k s]$ for $i \in \{1, 3\}$ and $k \in \{0, 2, 4, 6\}$. The group of order four generated by $[r, r^2 s]$ is normal in this subgroup, and the short exact sequence

$$\{e\} \longrightarrow \langle [r, r^2 s] \rangle \longrightarrow U \longrightarrow \langle [r^3, s] \rangle \longrightarrow \{e\}$$

splits, demonstrating that $U$ is isomorphic to $\mathbb{Z}_4 \rtimes \mathbb{Z}_2$. Conjugation of $[r, r^2 s]$ by $[r^3, s]$ yields $[r, r^2 s]^{-1}$, and therefore $U \cong D_4$. In conclusion, $\text{Aut}(\mathbb{Z}_8 \rtimes_\sigma \mathbb{Z}_2) \cong \mathbb{Z}_2 \times D_4$. $\square$

**Theorem 8.4.** $Aut(\mathbb{Z}_8 \rtimes_\tau \mathbb{Z}_2) \cong \mathbb{Z}_2 \times D_4$.

*Proof.* The element $[r^3, r^4 s]$ is central in $\text{Aut}(\mathbb{Z}_8 \rtimes_\tau \mathbb{Z}_2)$, and therefore $\langle [r^3, r^4 s] \rangle \cong \mathbb{Z}_2$ is normal in $\text{Aut}(\mathbb{Z}_8 \rtimes_\tau \mathbb{Z}_2)$. Letting $V$ denote the quotient $\text{Aut}(\mathbb{Z}_8 \rtimes_\tau \mathbb{Z}_2)/\langle [r^3, r^4 s] \rangle$, the short exact sequence

$$\{e\} \longrightarrow \langle [r^3, r^4 s] \rangle \longrightarrow \text{Aut}(\mathbb{Z}_8 \rtimes_\tau \mathbb{Z}_2) \longrightarrow V \longrightarrow \{e\}$$

splits, giving $\text{Aut}(\mathbb{Z}_8 \rtimes_\tau \mathbb{Z}_2) \cong \mathbb{Z}_2 \times V$.

The image of $V$ in $\text{Aut}(\mathbb{Z}_8 \rtimes_\tau \mathbb{Z}_2)$ through one possible section gives a subgroup containing elements of the form $[r^i s^j, s]$ for $i \in \{1, 3, 5, 7\}$ and $j \in \{0, 1\}$. The group of order four generated by $[r^3 s, s]$ is normal in this subgroup, and the short exact sequence

$$\{e\} \longrightarrow \langle [r^3 s, s] \rangle \longrightarrow V \longrightarrow \langle [r^3, s] \rangle \longrightarrow \{e\}$$

splits, demonstrating that $V$ is isomorphic to $\mathbb{Z}_4 \rtimes \mathbb{Z}_2$. Conjugation of $[r^3s, s]$ by $[r^3, s]$ yields $[r^3s, s]^{-1}$, and therefore $V \cong D_4$. In conclusion, $\text{Aut}(\mathbb{Z}_8 \rtimes_\tau \mathbb{Z}_2) \cong \mathbb{Z}_2 \times D_4$. □

## 8.4. The structure of $\text{Aut}(D_8)$.

**Theorem 8.5.** $\text{Aut}(D_8) \cong (\mathbb{Z}_2 \times D_4) \rtimes \mathbb{Z}_2$.

*Proof.* In $\text{Aut}(D_8)$ there is an index 2 subgroup $W$ containing elements of the form $[r^i, r^k s]$ for $i \in \{1, 3, 5, 7\}$ and $k \in \{0, 2, 4, 6\}$. $W$ is normal in $\text{Aut}(D_8)$ and the short exact sequence

$$\{e\} \longrightarrow W \longrightarrow \text{Aut}(D_8) \longrightarrow \langle [r^7, rs] \rangle \longrightarrow \{e\}$$

splits, giving $\text{Aut}(D_8) \cong W \rtimes \mathbb{Z}_2$. The action of $\langle [r^7, rs] \rangle$ on the members of $W$ will be described after the proof.

The automorphism $[r^5, s] \in \text{Aut}(D_8)$ generates a group of order 2 that is normal in $W$, and the sequence

$$\{e\} \longrightarrow \langle [r^5, s] \rangle \longrightarrow W \longrightarrow W/\langle [r^5, s] \rangle \longrightarrow \{e\}$$

splits. The image of the quotient $W/\langle [r^5, s] \rangle$ through one possible section gives a subgroup of $W$ containing elements of the form $[r^i, r^k s]$ for $i \in \{1, 3\}$ and $k \in \{0, 2, 4, 6\}$. The group of order four generated by $[r, r^2 s]$ is normal in this subgroup, and the short exact sequence

$$\{e\} \longrightarrow \langle [r, r^2 s] \rangle \longrightarrow W/\langle [r^5, s] \rangle \longrightarrow \langle [r^3, s] \rangle \longrightarrow \{e\}$$

splits, demonstrating that $W/\langle [r^5, s] \rangle$ is isomorphic to $\mathbb{Z}_4 \rtimes \mathbb{Z}_2$. Conjugation of $[r, r^2 s]$ by $[r^3, s]$ yields $[r, r^2 s]^{-1}$, and therefore $W/\langle [r^5, s] \rangle \cong D_4$.

Finally, $W \cong \mathbb{Z}_2 \times D_4$ and $\text{Aut}(D_8) \cong (\mathbb{Z}_2 \times D_4) \rtimes \mathbb{Z}_2$. □

The set $\{[r, r^2 s], [r^3, s], [r^5, s]\}$ generates $W$, and therefore the action of $\langle [r^7, rs] \rangle$ on $W$ in the product $\text{Aut}(D_8) \cong W \rtimes \langle [r^7, rs] \rangle$ can be seen via the image of each generator through the map "conjugation by $[r^7, rs]$":

$$[r^7, rs][r, r^2 s][r^7, rs] = [r, r^6 s]$$
$$[r^7, rs][r^3, s][r^7, rs] = [r^3, r^6 s]$$
$$[r^7, rs][r^5, s][r^7, rs] = [r^5, r^4 s]$$

By Theorem 7.2, $\text{Aut}(D_8) \cong \mathbb{Z}_8 \rtimes_\psi \text{Aut}(\mathbb{Z}_8)$ for $\psi$ isomorphic, and by the above Theorem 8.5 we have $\text{Aut}(D_8) \cong (\mathbb{Z}_2 \times D_4) \rtimes \mathbb{Z}_2$. We will show that these descriptions match by showing that the elements supplied in Theorem 8.5 have the structure $\mathbb{Z}_8 \rtimes_\psi \text{Aut}(\mathbb{Z}_8)$ given by Theorem 7.2.

The element $[r, rs]$ from $\text{Aut}(D_8)$ has order 8. The inverse of a given element $[r^i, r^k s]$ is $[r^i, r^{-ik} s]$, and conjugation of $[r, rs]$ by $[r^i, r^k s]$ gives

$$[r^i, r^k s][r, rs][r^i, r^{-ik} s] = [r, r^i s]$$

which demonstrates that the group $\mathbb{Z}_8$ generated by $[r, rs]$ is normal in $\text{Aut}(D_8)$. There are four cosets of $\langle [r, rs] \rangle$ in $\text{Aut}(D_8)$. They contain elements of the form $[r, r^k s]$, $[r^3, r^k s]$, $[r^5, r^k s]$ and $[r^7, r^k s]$, respectively. This quotient group is isomorphic to Klein's Group, and the short exact sequence $\{e\} \to \mathbb{Z}_8 \to \text{Aut}(D_8) \to \mathbb{Z}_2 \times \mathbb{Z}_2 \to \{e\}$ splits via a section that sends each coset of $\text{Aut}(D_8)/\langle [r, rs] \rangle$ to one of the elements $[r,s]$, $[r^i,s]$, $[r^5,s]$ and $[r^7,s]$. Klein's group is isomorphic to $\mathrm{Aut}(\mathbb{Z}_8)$, and conguagtion of $[r,rs]$ by $[r^i,s]$ gives $[r,r^i s]$, demonstrating (again) that $\mathrm{Aut}(D_8) \cong \mathbb{Z}_8 \rtimes_\psi \mathrm{Aut}(\mathbb{Z}_8)$ where $\psi$ is the identity isomorphism.

**Acknowledgments.** Many thanks to my mentor Bena Tshishiku for making this paper possible.